\documentclass[a4]{article}

\usepackage[francais]{babel}    % typo francaise
\usepackage[T1]{fontenc}        % accents codés dans la fonte
\usepackage[latin1]{inputenc}   % accents 8 bits dans le fichier
\usepackage{vmargin}            % régler la taille de la feuille

\usepackage{amsmath}    % les symboles les plus fréquents
\usepackage{amssymb}    % des symboles
\usepackage{amsfonts}   % des fontes, eg pour \mathbb

\setmarginsrb{3.5cm}{1cm}{3.5cm}{1.5cm}{14mm}{6mm}{0mm}{10mm}

\author{}
\title{$H^1$ n'a pas de base complètement inconditionnelle}
\begin{document}

\def\tens{\otimes}
\def\T{\mathbb T}
\def\N{\mathbb N}
\def\sk{\smallskip}
\renewcommand{\le}{\leqslant}
\renewcommand{\ge}{\geqslant}

\def\res #1#2{ {\leftskip=10mm \rightskip=10mm
\begin{center}{\bf #1}\end{center} \vskip -2mm
{\small \indent #2} \par}}

\newcommand{\theo}[2]{
\sk
\noindent{\bf #1} {\it #2}
\sk
}
 
\maketitle

\res{Résumé}{\indent On montre que l'espace de Hardy $H^1$ n'admet pas de 
décomposition inconditionnelle de rang fini de l'identité
en tant qu'espace d'opérateurs. Il en résulte
qu'il n'existe pas de base complètement inconditionnelle dans $H^1$.}

\sk

\begin{center}
{\bf $H^1$ does not have any completely unconditional basis}
\end{center}
\res{English summary}{\indent Let $H^1$ be the classical Hardy space 
of analytic functions on
the unit disc. We show that this space does not admit any finite rank 
completely unconditional decomposition of the identity, 
as a consequence 
it fails to have a completely unconditional basis.
}  

\sk
\sk

L'objet de cette note est de montrer que $H^1$ n'admet pas de base 
complètement inconditionnelle et plus
généralement pas de décomposition complètement 
inconditionnelle de rang fini. On dit qu'un espace de Banach $X$ admet
une décompostion inconditionnelle de rang fini s'il existe une suite 
d'endomorphismes de X de rang fini $(P_n)$ tels que
$$ \sup_{\epsilon_n=\pm 1} \sup_N \left\|\sum_{n=0}^N \epsilon_n 
P_n\right\|_{X\to X}
<\infty$$ 
et pour tout $x\in X$, $x=\sum_{n=0}^\infty P_n(x)$ où la
série converge inconditionnellement. Si $X$ est un 
espace d'opérateurs, une telle décompostion est dite complètement
inconditionnelle si en outre
 $$ \sup_{\epsilon_n=\pm 1} \sup_N \left\|\sum_{n=0}^N 
\epsilon_n P_n\right\|_{cb(X, X)}
<\infty.$$ 
 $H^1$ désigne l'espace de Hardy du tore, il s'agit du sous-espace de
$L^1(\T)$, formé des fonctions intégrables à valeurs complexes sur le cercle 
unité du plan complexe
ayant tous leurs coefficients de Fourier d'indice négatifs nuls.
 Cet espace de Banach admet d'autres descriptions en analyse réelle, 
par exemple il est constitué de l'ensemble des fonctions du type
$f+i H(f)$, où $f$ est une fonction réelle intégrable sur $\T$ ayant
une transformée de Hilbert $H(f)$ également intégrable. L'autre 
caractérisation importante essentiellement 
dûe à Fefferman, établit que les parties 
réelles des fonctions dans $H^1$ se décomposent en atomes 
(voir \cite{Ga}).    
 
 Le premier résultat de théorie locale pour $H^1$ remonte à Stein 
(\cite{St}), les multiplicateurs de la Vallée-Poussin 
(convolution avec les noyaux $W_n$, $n\ge 0$ du même nom) forment une 
décomposition inconditionnelle de l'identité, c'est à dire qu'il existe
$K\ge 0$ telle que pour tout choix de signes $\epsilon_n$ et tout 
$N\ge 0$ :
$$\left\|f\ast (\sum_{n=0}^N \epsilon_n W_n)\right\|_{H^1} \le K
 \,\left\|f\right\|_{H^1}$$
et $f\ast (\sum_{n=0}^N W_n) $ tend vers $f$ dans $H^1$ lorsque $N$ tend
vers l'infini. L'étape suivante était de savoir si cet espace de Banach 
admet une base inconditionnelle, Maurey dans 
\cite{Ma} y a répondu par l'affirmative mais d'une façon assez 
indirecte.
 Ensuite, Wojtaszczyk (\cite W), s'inspirant de travaux de Carleson, a
explicité une telle base ; le système de Franklin s'avère être une base
inconditionnelle de l'espace atomique $H^1(\mathbb R)$. 
Depuis, la théorie
des ondelettes a permis de mieux comprendre la situation, voir à 
ce sujet le livre 
d'Yves Meyer \cite{Me}, chapitres V et VI.

 L'espace $H^1$ est muni d'une structure d'espace d'opérateurs 
au sens de la théorie développée par Blecher-Paulsen et Effros-Ruan
({\it cf} \cite{ER}), 
en d'autre termes $H^1$ peut être réalisé comme un sous-espace 
de l'espace $B(\mathcal H)$
des opérateurs bornés d'un Hilbert $\mathcal H$
dans lui-même (avec $\mathcal H=\ell_2$ 
car $H^1$ est séparable). On note $S^1$ ($S^1_d$) l'espace 
des opérateurs à trace
de $\ell_2$ ($\ell_2^d$) muni de sa norme usuelle,
et $H^1(S^1)$ l'ensemble des fonctions de $\T$
à valeurs dans $S^1$ ayant des coefficients de Fourier 
négatifs nuls, muni de la norme
induite par $L^1(S^1)$.
La question de l'existence de bases 
inconditionnelles admet un analogue dans cette théorie. En utilisant
les résultats de \cite{Pi}, elle se retraduit de la manière suivante :
peut-on trouver une base de $H^1$, $\{\psi_1,\dots,\psi_n,\dots\}$ et
 $K\ge 0$ 
telles que pour tout $N$ entier et toute famille $\{x_i\}$ dans $S_1$
 et 
tout choix de signes $\epsilon_n$:
$$\left\| \sum_{n=0}^N \epsilon_n x_n\tens \psi_n\right\|_{H^1(S^1)}\le
K \, \left\|\sum_{n=0}^N x_n\tens \psi_n\right\|_{H^1(S^1)} \,\,?$$ 
 Pour infirmer cela, on démontre qu'il n'y a pas d'analogue du résultat 
de Stein.

\theo{Théorème }{L'espace $H^1(\T)$ 
n'admet pas de décomposition complètement inconditionnelle 
de rang fini de l'identité.}
 
 Cela signifie qu'on ne peut pas trouver de décomposition 
inconditionnelle $(P_n)$ de rang fini de $H^1$,
 telle que la série $\sum P_n\tens Id_{S^1}$
soit inconditionnelle dans $\mathcal B(H^1(S^1),H^1(S^1))$.
  
Plus précisément on montre que si $(P_n)$ est une telle décomposition,
alors  
on peut construire une suite $a_n$ bornée en valeur absolue par 1, 
telle que
$$\left\|\sum a_n P_n\tens Id\right\|_{H^1(S^1_d)\to H^1(S^1_d)}\ge K \ln d.$$

 On a besoin du fait élémentaire suivant :
 
\theo{Lemme}{Soit $\sum f_k$
 une série inconditionnelle dans $L_1(\T)$,
alors les séries $\sum \epsilon_k f_k$ sont uniformément de Cauchy
en $\epsilon_k\in [-1,1]$.}

{\it Démonstration du théorème :}
 On adopte la notation, $e_n(t)=e^{2\pi int}$. La preuve 
s'inspire d'un des résultats de \cite N. L'idée est la suivante :
à partir d'une décomposition inconditionnelle de l'identité $(P_n)$, 
on construit des sortes de 
multiplicateurs qui permettent ensuite de transférer la projection 
triangulaire de $S^1_d$ comme endomorphisme de $H^1(S^1_d)$.
\\ Soit $\eta>0$ fixé, on construit par récurrence 
\\ \indent $(\phi_n)$ une suite d'applications de rang fini de la 
forme $\sum a_k P_k$ avec $a_k\in\{0,1\}$,
\\ \indent $(\alpha_n),(\beta_n)$ des suites croissantes d'entiers,
\\ \indent $(\epsilon_n)$ une suite de réels croissante majorée par
$\eta$, 
\\ telles que,
$$ \forall \, i,j\le n \left\{\begin{array}{cc}
i) \,\,j\le i & \left\|\phi_n(e_{\alpha_i+\beta_j})\right\|\le \epsilon_n\\
ii)\,\,i< j   & 
\left\|\phi_n(e_{\alpha_i+\beta_j})-e_{\alpha_i+\beta_j}\right\|\le \epsilon_n.
\end{array}\right. $$

On initialise, la récurrence par $\alpha_1=\beta_1=0$,
 $\epsilon_1=\epsilon$
un réel quelconque (petit). On choisit $\phi_1=0$, $(i)$ et $(ii)$ 
sont vérifiées.

\sk

On suppose toutes les suites construites jusqu'au rang $n$.
\\ On fixe $\delta>0$.
D'après le lemme, on peut trouver un entier $K$, tel que
$$\forall\, i,j\le n \quad \forall\, l\ge k \ge K \quad 
\forall \, a_k\in \{0,1\},
\quad \left\|\sum_{t=k}^l a_t P_t(e_{\alpha_i+\beta_j})\right\|\le \delta.$$   
On peut supposer que si $\phi_n=\sum_{k=1}^A a_k P_k$, alors $K>A$.

\sk 

On détermine $\beta_{n+1}$ ; 
\\ \indent puisque les $P_k$ sont de rang fini, on a $\lim_{\beta\to \infty}
P_k(e_\beta)=0$ dans $H^1$ car $e_\beta$ tend faiblement vers 0.
On choisit $\beta_{n+1}$ suffisament grand pour que pour tout 
$\beta \ge \beta_{n+1}$ et tout choix $(a_k)_1^K\in \{0,1\}^K$ 
(il y en a un nombre fini)
$$\left\|\sum_{k=1}^K a_k P_k(e_\beta)\right\|\,< \,\delta.$$
On détermine $\phi_{n+1}$ ; 
\\ \indent puisque $\sum P_k$ tend simplement vers l'identité de $H^1$, 
on peut trouver un entier $N>K$ de sorte que pour tout $i\le n$ 
$$\left\| \sum_{k=1}^N P_k(e_{\alpha_i+\beta_{n+1}}) \,-\,e_{\alpha_i+\beta_{n+1}}
\right\|<\delta.$$
On pose $\phi_{n+1}=\phi_n+\sum_{k=K}^N P_k$. Le choix de $K$ assure 
que $\phi_{n+1}$ est de la forme 
$\sum a_k P_k$ avec $a_k\in\{0,1\}$.
\\ \indent On vérifie $(ii)$ pour $i<n+1$
$$\left.\begin{array}{lcl}
\left\|\phi_{n+1}(e_{\alpha_i+\beta_{n+1}})-e_{\alpha_i+\beta_{n+1}}\right\|
& \le & \left\|\phi_{n}(e_{\alpha_i+\beta_{n+1}})\right\| +
 \left\|\sum_{k=1}^{K-1} P_k
(e_{\alpha_i+\beta_{n+1}})\right\|\\
& & +\left\|\sum_{k=1}^N P_k
(e_{\alpha_i+\beta_{n+1}})-e_{\alpha_i+\beta_{n+1}}\right\|\\
 & \le & \delta+\delta +\delta \\
 &\le &3 \delta.
\end{array}\right.$$
On vérifie $(ii)$ pour $i<j\le n$
$$\left.\begin{array}{lcl}
\left\|\phi_{n+1}(e_{\alpha_i+\beta_{j}})-e_{\alpha_i+\beta_{j}}\right\|
& \le &\left\|\phi_{n}(e_{\alpha_i+\beta_{j}})-e_{\alpha_i+\beta_{j}}\right\|
 +\left\|\sum_{k=K}^N P_k(e_{\alpha_i+\beta_{j}})\right\|
\\ & \le & \epsilon_n+\delta.
\end{array}\right. $$
Il reste à fixer $\alpha_{n+1}$ et vérifier $(i)$.
\\ \indent Puisque $\phi_{n+1}$ est de rang fini, 
$\lim_{\alpha\to\infty}
\phi_{n+1}(e_\alpha)=0$ dans $H^1$, on peut donc choisir $\alpha_{n+1}$
de sorte que $(ii)$ soit vérifiée pour $i=n+1$
$$\forall \, j\le n+1 \qquad 
\left\|\phi_{n+1}(e_{\alpha_{n+1}+\beta_{j}})\right\|
\le \delta.$$
 Dans le cas où $j\le i\le n$, on a 
$$\left.\begin{array}{lcl}
\left\|\phi_{n+1}(e_{\alpha_{i}+\beta_{j}})\right\|& \le& 
\left\|\phi_{n}(e_{\alpha_{i}+\beta_{j}})\right\|+\left\|\sum_{k=K}^N 
P_k(e_{\alpha_i+\beta_{j}})\right\|\\
 & \le & \epsilon_n + \delta.
\end{array}\right. $$
 Finalement, avec $\epsilon_{n+1}=\max (3\delta, \epsilon_n+\delta)$,
l'hypothèse de récurrence est vérifiée pour un choix
convenable de $\delta$, et on peut en outre supposer 
que $\epsilon_{n+1}\le (1+2^{-(n+1)})\epsilon_n$.

\sk 

 En résumé, il est possible de construire des suites avec les propriétés
$(i)$ et $(ii)$, avec $\epsilon_n\le \epsilon\prod_{k=1}^n (1+2^{-k})<\eta$
pour $\epsilon$ suffisamment petit.

\sk

On utilise cette construction pour minorer 
$$C=\sup_{a_k\in \{0,1\} }
\left\|\sum a_k P_k\otimes Id_{S^1_d}\right\|_{H^1(S^1_d)
\to H^1(S^1_d)} \qquad \qquad (*) $$
En particulier, les $\left\|\phi_n\otimes Id_{S^1_d}\right\|$
avec les $\phi_n$ 
construits précédemment sont bornées par $C$.
Soit $X=(x_{i,j})$ une matrice de taille $d$, l'élément de $H^1(S^1_d)$ 
défini par 
$$Z=diag(e_{\alpha_i})Xdiag(e_{\beta_i}),$$
où $diag(d_i)$ est la matrice diagonale de taille $d$ ayant pour entrées les
$d_i$, vérifie $$\left\|Z\right\|_{H^1(S^1_d)}=\left\|X\right\|_{S^1_d}.$$
 \\ On note $T$ la projection triangulaire supérieure dans $S^1_d$ ;
il s'agit de l'application linéaire de $S^1_d$ dans lui-même qui
à une matrice $X=(x_{i,j})_{1\le i, j\le d}$ associe la 
matrice 
$$T(X)=(x_{i,j}1_{\{j>i\}}).$$
 Les 
propriétées $(i)$ et $(ii)$ impliquent
$$\left\|\phi_n\otimes 
Id_{S^1_d}(Z)-Id_{H^1}\otimes T(Z)\right\|_{H^1(S^1_d)}\le
\epsilon_n d^2 \left\|X\right\|_{S^1_d}.$$
D'autre part, $Id_{H^1}\otimes T(Z)=diag(e_{\alpha_i})T(X)diag(e_{\beta_i})$
et donc $\left\|Id_{H^1}\otimes T(Z)\right\|_{H^1(S^1_d)}=
\left\|T(X)\right\|_{S^1_d}$, d'où
$$  \left\|T(X)\right\|_{S^1_d}\le \epsilon_n d^2 \left\|X\right\|_{S^1_d} 
+ C \left\|X\right\|_{S^1_d}.$$
$\epsilon_n$ pouvant être arbitrairement petit, on en déduit que 
$C$ domine la norme de la projection triangulaire dans $S^1_ d$ 
qui est de l'ordre
de $\ln d$ d'après \cite {Kp}, ce qui termine la preuve.

\theo{Corollaire}{
$H^1$ n'admet pas de base complètement inconditionnelle.}
 
Si tel était le cas, les projections sur les vecteurs de base formeraient
une décompostion complètement inconditionnelle de l'identité.

 Il découle de \cite{N} ou de \cite {Pe} 
que l'estimation en $\ln d$ pour la contante
$(*)$ est optimale et atteinte pour la décomposition de Stein.

\noindent \'Eric Ricard
\\Universit{\'e} Paris VI
\\{\'E}quipe d'analyse
\\4 place Jussieu
\\75252 Paris Cedex 05 
\\France


\begin{thebibliography}{99}

\bibitem{ER} Effros E. et Ruan Z.J. {\it Operator spaces.} Oxford Univ.
Press. À paraître. 

\bibitem{Ga} Garnett, J. {\it Bounded analytic functions}, Academic Press, 
New York, 1981.

\bibitem{Kp} Kwapie\'n, S. et Pe\l czy\'nski A.  
The main triangle projection in matrix spaces and its
applications. {\it Studia Math}. 34 (1970) 43-68.

\bibitem{Ma} Maurey, B. Isomorphismes entre espace $H^1$,
{\it Acta math.} 145 (1980) p 79-120.

\bibitem{Me} Meyer Y. {\it Ondelettes} Hermann 1990. 

\bibitem{N} Nazarov F., Pisier G., Treil S., Volberg A. Sharp 
estimates in vector Carleson imbedding theorem and for vector
paraproducts. À paraître.  

\bibitem{Pe} Petermichl S. Dyadic shift and a logarithmic estimate
for Hankel operator with matrix symbol,
{\it C. R. Acad. Sci. Paris} 2000, À paraître.

\bibitem{Pi} Pisier, G. {\it Non-commutative vector valued $L_p$-spaces and
completely $p$-summing maps}, Astérisque 247 (1998).

\bibitem{St} Stein, E. Multiplicateurs et fonctions de Littlewood-Paley, 
{\it C. R. Acad. Sci. Paris Série A-B}, 263 (1966) A 716-719.

\bibitem{W} Wojtaszczyk, P. The Franklin system is an unconditional basis
in $H^1$, {\it Arkiv för Mat.}, 20 (1982), 293-300.

\end{thebibliography}
\end{document}